
%
%
%
%
%
\magnification=1200
\baselineskip=16.5pt plus .2pt  
\parskip=4.5pt plus1.5pt minus 1pt  
\def\cite{\bf}
\def\theorem{\vskip\baselineskip\bf}
\def\demo{\vskip\baselineskip\par{\bf Proof:}}
\def\qed{{\bf Q.E.D.}}
\centerline {\bf A Non-Spectral Dense Banach Subalgebra
of the Irrational Rotation Algebra}
\vskip\baselineskip
\centerline {\bf Larry B. Schweitzer}
\vskip\baselineskip
\centerline {{\it Current address:}\ \  Department of
Mathematics, University of Victoria}
\centerline{  Victoria, B.C. Canada V8W 3P4}
\centerline{{\it Email Address:}\ \ lschweit@sol.uvic.ca}
\vskip\baselineskip
\centerline {{\it 1991 Mathematics Subject Classification.}
Primary: 46H99, Secondary: 46H35, 46H25, 46L99}
\vskip\baselineskip
\centerline{\bf Abstract}
\par
We  give an example of a dense, simple, unital
Banach subalgebra $A$
of the irrational rotation C*-algebra $B$, such that $A$ is not
a spectral subalgebra of $B$.
This  answers a question posed in  T.W. Palmer's paper
{\cite{[1]}}.
\vskip\baselineskip
%
%
%
%
\par
If $A$ is a subalgebra of an algebra $B$ (both algebras over
the complex
numbers), we say that $A$ is a
{\it spectral subalgebra} of $B$
if the quasi-invertible
elements of $A$ are precisely the quasi-invertible elements of
$B$ which lie in $A$.   In the language of {\cite{[3]}},
this is equivalent to saying that $A$ is a {\it spectral
invariant
subalgebra} of $B$.
\par
There are many known examples of dense unital Banach
subalgebras of
C*-algebras which are not spectral.  For example, see
Example 3.1 of {\cite{[3]}}.   The example we give here
is of interest
because the Banach algebra is simple, and thus answers
Question 5.12 of {\cite{[1]}} in the negative.
\par
Recall that the irrational rotation  algebra associated
with
an irrational real number $\theta$ is the C*-crossed
product of the integers $\bf Z$ with the commutative
C*-algebra of continuous functions on the circle
$C(\bf T )$, where $n \in \bf Z$ acts via
$ \alpha_n(\varphi)(z) = \varphi(z-n\theta)$, for
$\varphi \in C(\bf T)$
and $z\in \bf T$.
Let $B= C^{*}(\bf Z,  C(\bf T), \theta)$ denote
this crossed product.
\par
Let $A$ be the set of functions F from $\bf Z$
to $C(\bf T)$ which satisfy the
integrability condition
$$ \parallel F \parallel_{A} =
\sum_{n \in \bf Z} e^{|n|} \parallel F(n)
\parallel_{\infty} <\infty, $$
where $\parallel \quad \parallel_{\infty}$ denotes the
sup norm on $C(\bf T)$.
Then $A$ is complete for the norm
$\parallel \quad \parallel_{A}$ and is a Banach
algebra.  The algebra
$A$ is contained in $L^{1}(\bf Z, C(\bf T))$ with dense
and continuous
inclusion, and hence contained in $B$
with dense and continuous inclusion.
Recall that the multiplication (in both $A$ and $B$)
is given by
$$ F*G(n, z) = \sum_{m \in \bf Z} F(n, z) G(n-m, z-m\theta),
\qquad
F, G \in A, \quad n \in \bf Z,\quad z\in \bf T.  $$
Let $u_n =\delta_{n} \otimes 1 \in A$ denote the delta
function at
$n\in \bf Z$
tensored with the identity in $C(\bf T)$.
Then $u_0 $ is the unit in both
$A$ and $B$.
{\theorem{Theorem 1}}{\it{ The Banach algebra $A$ is simple.}}
{\demo}
We imitate the argument of  {\cite{[2]}}.
Define a continuous linear map $P\colon A \rightarrow
C(\bf T)\subseteq A$ by $P(F) = F(0)$.   Note that
$\parallel P(F) \parallel_{A}
\leq \parallel F\parallel_{A}$ for $F \in A$.
Let $J$ be a closed two-sided ideal in $A$, which is not
equal to $A$.  Since
$\bf Z$ acts ergodically on $\bf T$, we know
that $C(\bf T)$ has no nontrivial
closed $\bf Z$-invariant
ideals. Hence $J\cap C(\bf T) = 0$.
\par
We show that $P(J)=0$.
It suffices to show that $P(J) \subseteq J$.
Let $\epsilon >0$ and $F \in A$. Let $N$ be a sufficiently
large
integer for which
$$ \sum_{|n|>N} e^{|n|} \parallel F(n) \parallel_{\infty}
< \epsilon. $$
Define $F_{1} \in A$ by $F_{1}(n) = 0 $ if $|n|>N$, and
$F_{1}(n)= F(n)$
for $|n| \leq N$.
By the proof of Lemma 6 of {\cite{[2]}}, there exists
unimodular functions
$\theta_{1},\dots \theta_{M} \in C(\bf T)$ such that
$$ P(F_{1}) = {1\over{M}} \sum_{n = 1}^{M}
\theta_{n}^{*} F_{1} \theta_{n}.$$
(Here unimodular means that $|\theta_{i}(z)|=1$ for
each $z \in \bf T$
and $i = 1, \dots M$.)
Hence
$$\parallel P(F) - {1\over{M}} \sum_{n=1}^{M}
\theta_{n}^{*} F \theta_{n}\parallel_{A} \,\,
\leq\,\, \parallel P(F-F_{1}) \parallel_{A} \,\,+ \,\,
\parallel F-F_{1} \parallel_{A} \,\,
< \,\,2\epsilon. \eqno (*) $$
Now if $F\in J$, (*) shows that $P(F)$ can be
approximated arbitrarily
closely by elements of $J$.
Since $J$ is closed, this shows that $P(F) \in J$.
Hence $P(J) \subseteq J$ and $P(J)=0$.
\par
 If $P(Fu_n )=0$ for all $n$, then $F(n)= 0$
for all $n$ and so $F = 0$. Since $J$ is a
two-sided ideal
and $P(J)=0$, we have $P(Ju_n)=0$ for all $n$.  Hence $J=0$
and $A$ is
simple.
\qed
{\theorem{Theorem 2}} {\it The Banach algebra $A$ is not a
spectral
subalgebra of $B$. }
{\demo}  We construct an algebraically irreducible
$A$-module which is not
contained in any *-representation of $B$ on a
Hilbert space.
By  Corollary 1.5 of {\cite{[3]}}, it will follow that
$A$ is not
a spectral subalgebra of $B$.
\par
Let $E$ be the Banach $A$-module $C(\bf T)$ with sup norm,
and with
(continuous) action of $A$ given by
$$
(F\varphi)(z) = \sum_{n} F(n, z) e^{n} \varphi(z - n \theta),
\qquad \varphi \in E,\quad F\in A,\quad z\in \bf T.
$$
We show that $E$ is in fact algebraically irreducible.
Let $\varphi \in E$
be not identically equal to zero.  Since the complex
conjugate of $\varphi$
is in $A$, the algebraic span $A \varphi$ contains
$|\varphi|^{2}$,
which we denote by $\psi$.
Note $u_n \psi(z) =
e^{n} \psi(z-n\theta)$.  Since $\theta$ is irrational
and $\bf T$ is
compact, there exists finitely many $n_{1}, \dots n_{k}
\in \bf Z$
such that the sum of $u_{n_i}\psi$ from  $i=1$
to $k$ never vanishes on $\bf T$.  If $\chi$ is this sum,
then $1/\chi$
is in $C(\bf T) \subseteq A$
so $1 \in A\varphi$ and hence $E = A \varphi$.   This
proves that $E$ is algebraically irreducible.
\par
It remains to show that no *-representation of $B$ on a
Hilbert space contains $E$.
But the  action of $\bf Z$ on $1\in E$ is given  by
$u_n 1 = e^n 1$.  Clearly the Hilbert space could not
have a unitary,
or even isometric, action of $\bf Z$.
\qed
\vskip\baselineskip
\centerline{\bf References}
\vskip\baselineskip
\frenchspacing
\item{[1]}  T.W. Palmer, {\it Spectral Algebras},
Rocky Mtn. Jour. Math.  22(1) (1992), 293-328.
\frenchspacing
\item{[2]} S.C. Power, {\it Simplicity of C*-algebras of
Minimal Dynamical Systems}, J. London Math. Soc. (2) 18
(1978), 534-538.
\frenchspacing
\item{[3]} L.B. Schweitzer,
{\it A Short Proof that $M_n(A)$ is local if
A is local and Fr\'echet},
Intl. Jour. Math. 3 (4) (1992), 581-589.
\vskip\baselineskip
\end